\documentclass[11pt,reqno]{amsart} 
 \usepackage{amsmath,amssymb,amsthm}
 \usepackage{graphicx}
 \usepackage{cite}
 \usepackage[mathscr]{eucal}

 \theoremstyle{plain}
 \newtheorem{theorem}{Theorem}[section]  
   
 \newtheorem{lemma}[theorem]{Lemma}  
 
 \newtheorem*{theorem*}{Theorem}

 \theoremstyle{plain}

 \newtheoremstyle{citing}
   {3pt}
   {3pt}
   {\itshape}
   {}
   {\bfseries}
   {.}
   {.5em}
   {\thmnote{#3}}

 \theoremstyle{citing}

\numberwithin{equation}{section}
\allowdisplaybreaks[1]

\newlength{\intwidth}
\makeatletter
\DeclareRobustCommand{\cpvint}[2]
    {\mathop{%
       \text{%
         \settowidth{\intwidth}{%
           \ifx\ilimits@\displaylimits
             $\int_{#1}^{#2}$%
           \else
             $\int$%
           \fi}%
         \makebox[0pt][l]{\makebox[\intwidth]{$\text{C}$}}%
         $\int_{#1}^{#2}$}}}
\makeatother

\makeatletter
\DeclareRobustCommand{\cpvintsmall}[2]
    {\mathop{%
       \text{%
         \settowidth{\intwidth}{%
           \ifx\ilimits@\displaylimits
             $\int_{#1}^{#2}$%
           \else
             $\int$%
           \fi}%
         \makebox[0pt][l]{\makebox[\intwidth]{$\text{{\tiny C}}$}}%
         $\int_{#1}^{#2}$}}}
\makeatother

\newcommand{\where}{:\:}

\newcommand{\nz}{{\mathbb N}}

\newcommand{\rz}{{\mathbb R}}

\renewcommand{\phi}{\varphi} 
\newcommand{\eval}{\vert}

\begin{document}
 
\title[Planar Plateau Problem] {Multiple Solutions to the Planar Plateau Problem}
\author{Matthias Schneider}
\address{Ruprecht-Karls-Universit\"at\\
         Im Neuenheimer Feld 288\\
         69120 Heidelberg, Germany\\}
\email{mschneid@mathi.uni-heidelberg.de} 
\date{February 16, 2010}  

\begin{abstract}
We give existence and nonuniqueness results for simple planar 
curves with prescribed geodesic curvature.
\end{abstract}
\maketitle

\section{Introduction}
\label{sec:introduction}
We are interested in the planar Plateau problem:
Given two points $p_1$ and $p_2$ in the plane and a smooth
function $k:\rz^2\times [0,1]\to \rz$, find an immersed curve $\gamma \in C^2([0,1],\rz^2)$,
such that $\gamma(0)=p_1$, $\gamma(1)=p_2$, and for every $t\in [0,1]$
the (signed) geodesic curvature $k_\gamma(t)$ of $\gamma$ at $t$, 
\begin{align*}
k_\gamma(t) := |\dot\gamma(t)|^{-3}\big\langle \ddot\gamma(t) ,J\dot\gamma(t)\big\rangle,  
\end{align*}
is given by $k(\gamma(t),t)$, where $J$ denotes the rotation by $\pi/2$.
We choose the orientation, such that the circle
of radius $r$ with counterclockwise parameterization has positive curvature $r^{-1}$.\\
Without loss of generality after a rotation and a translation we may
assume that $p_1= (a,0)$ and $p_2=(-a,0)$ for some $a>0$.
Then the planar Plateau problem is equivalent to the following
ordinary differential equation
\begin{align}
\label{eq:1}
\ddot \gamma = |\dot\gamma| k(\gamma(t),t) J(\dot\gamma),\\
\notag \gamma(0)=(a,0),\, \gamma(1)=(-a,0),  
\end{align}
If the function $k\equiv k_0$ is constant, by elementary geometry, 
the planar Plateau problem is only solvable for $|k_0|\le a^{-1}$;
the solutions in this case 
are given by subarcs connecting $(a,0)$ and $(-a,0)$ of $n$-fold iterates of
a circle of radius $|k_0|$ with clockwise or counterclockwise parameterization depending
on the sign of $k_0$.  
If the analysis is restricted to simple solutions, then there
are $2$ solutions if $|k_0|<a^{-1}$, the {\em small} and the {\em large}  
solution corresponding to the subarcs subtending an angle strictly smaller or strictly larger than $\pi$. 
If $k_0=\pm a^{-1}$
then the unique simple solution is given by the half circle lying above
or below the $x$-axis depending on the sign of $k_0$. We will be mainly
interested in the case when $k$ is a positive function.\\
If the prescribed curvature function is independent of the variable $t$,
then the planar Plateau problem is 'geometric', in the sense that the set
of solutions is invariant under reparameterizations. If in this case
the function $k$ satisfies
$\|k\|_\infty < a^{-1}$, then 
from \cite{MR2010409} there exists a stable solution $\gamma_s$ to $(P)$.
We refer to $\gamma_s$ as a {\em small} solution. In the higher dimensional case
and in the context of $H$-surfaces analogous results are given in \cite{MR737190,MR0250208}.
For closed curves with prescribed curvature we refer to 
\cite{MR2319450,MR2254070,MR2287704,sissa082010,arXiv:0808.4038}.\\ 
Concerning the existence of a second, {\em large} solution for non-constant functions
$k$ there are only perturbative results, i.e. the function $k$ is assumed to be
close to a constant $k_0$, see \cite{MR2010409}. Concerning the existence of a large $H$-surface
we refer to \cite{MR733715,MR823116,MR832287}, if $H$ is constant, and
to \cite{MR1262929,MR2221202,MR1039366,MR1305276,MR1306675,MR1198307} for non-constant functions $H$.\\
We give existence criteria
for a large solution, that are non-perturbative.  
\begin{theorem}
\label{intro:main_theorem}
Let $a>0$ and $k \in C(\rz^2\times [0,1],\rz)$ be given, such that
\begin{align*}
0<\inf_{\rz^2\times [0,1]} k \le \sup_{\rz^2\times [0,1]} k <a^{-1},
\end{align*}
then there is a simple curve that solves \eqref{eq:1}. If, moreover, 
\begin{align}
\label{eq:cond_pinch}
\frac
{\sup_{(x,t) \in \rz^2\times [0,1]} k(x,t)}
{\sup_{(x,t) \in \rz^2 \times [0,1]} k(x,t)a+1}
< 
\inf_{(x,t) \in \rz^2\times [0,1]} k(x,t),    
\end{align}
then equation \eqref{eq:1} possesses at least two simple solutions.
\end{theorem}
To illustrate the pinching condition \eqref{eq:cond_pinch}
we note that the assumptions of Theorem \ref{intro:main_theorem} are satisfied, if
\begin{align*}
\frac12 a^{-1}< \inf_{\rz^2\times [0,1]} k \text { and }
\sup_{\rz^2\times [0,1]} k<a^{-1}.  
\end{align*}
The small solution is found in the set
\begin{align*}
M_{small} &:=\big\{\gamma \in C^2([0,1],\rz^2) \where
\gamma(0)=(a,0),\, \gamma(1)=(-a,0),\\  
&\qquad \gamma \oplus [-a,a] \text{ is simple, }
|\dot\gamma(0)|^{-1}\dot\gamma(0) \in
\{e^{i\theta}\where \pi/2<\theta<\pi\},\\
&\qquad |\dot\gamma(1)|^{-1}\dot\gamma(1) \in \{e^{i\theta}\where \pi<\theta<3\pi/2\}\big\},  
\end{align*}
whereas the large solution belongs to
\begin{align*}
M_{large} &:=\bigg\{\gamma \in C^2([0,1],\rz^2) \where
\gamma(0)=(a,0),\, \gamma(1)=(-a,0),\\ 
&\qquad \gamma \oplus [-a,a] \text{ is simple, }
\frac{\dot\gamma(0)}{|\dot\gamma(0)|} \in
\{e^{i\theta}\where -\pi/2<\theta<\pi\},\\
&\qquad \frac{\dot\gamma(1)}{|\dot\gamma(1)|}  \in
\{e^{i\theta}\where \pi<\theta<5\pi/2\},
\text{ and}\\ 
&\qquad \bigg(\frac{\dot\gamma(0)}{|\dot\gamma(0)|} \in
\{e^{i\theta}\where -\pi/2<\theta<\pi/2\} 
\text{ or }\\
&\qquad \frac{\dot\gamma(1)}{|\dot\gamma(1)|}  \in
\{e^{i\theta}\where 3/2\pi<\theta<5\pi/2\}\bigg)
\bigg\},
\end{align*}
where we define for a curve $\gamma \in C^0([0,L],\rz^2)$ connecting $(a,0)$ and $(-a,0)$
the closed curve $\gamma\oplus[-a,a] \in C([0,L+2a],\rz^2)$ by
\begin{align*}
\gamma\oplus[-a,a](t) :=  
\begin{cases}
\gamma(t) & 0\le t\le L\\
(-a+t-L,0) &L\le t\le L+2a.  
\end{cases}
\end{align*}
\begin{figure}[ht]
\resizebox{11cm}{!}{\includegraphics{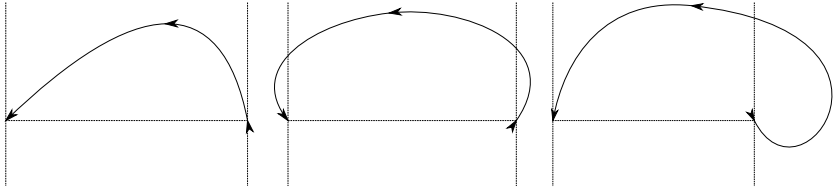}}
\caption{Examples of a small solution and two large solutions}
\label{fig:soln}
\end{figure}
The existence result is proved by using the Leray-Schauder degree
and suitable apriori estimates, i.e. we show that the degree of \eqref{eq:1}
with respect to $M_{small}$ equals $1$ and is given by $-1$, when computed in the set $M_{large}$.
The existence result then follows, since a non vanishing degree gives rise to a solution.     
The degree approach is interesting in itself and yields the flexibility
to deal with functions $k$ that depend on $x$ and $t$, for instance 
if $k$ does only depend on $t$, then the existence result shows
that in contrast to the four vertex theorem for simple closed 
curves of prescribed curvature (see \cite{MR2137880,MR2285124}) there is no additional condition on $k$
besides the $L^\infty$-bound for the corresponding boundary value problem.       
Moreover, the degree argument gives
the perspective to be applied to the higher dimensional case as well, e.g.
to surfaces in $\rz^3$ with prescribed mean curvature. 

\section{Apriori estimates}
\label{sec:apriori}

\begin{lemma}
\label{l:graph_small}  
Let $\gamma \in C^2([0,L],\rz^2)$ be a unit speed curve with positive
geodesic curvature connecting
$(a,0)$ and $(-a,0)$, such that the closed curve
$\gamma\oplus[-a,a] \in C([0,L+2a],\rz^2)$
is simple.\\
If 
\begin{align*}
\dot\gamma(0) &= e^{i\theta_0} \text{ for some } \pi/2\le \theta_0<\pi \text{ and }\\
\dot\gamma(L) &= e^{i\theta_L} \text{ for some } \pi<\theta_L\le \frac32\pi,  
\end{align*}
then $\gamma$ is a graph over the $x_1$-axis and there is a strictly decreasing
$C^2$-function
$\theta:[0,L] \to [\theta_L,\theta_0]$ such that
\begin{align*}
\dot\gamma(t)=e^{i\theta(t)}.
\end{align*}
If 
\begin{align*}
\dot\gamma(0)&= e^{i\theta_0} \text{ for some } -1/2\pi\le \theta_0 < \pi  \text{ and }\\
\dot\gamma(L) &= e^{i\theta_L} \text{ for some } \pi<\theta_L\le \frac52\pi,  
\end{align*}
then there are a strictly increasing
$C^2$-function
$\theta:[0,L] \to [\theta_0,\theta_L]$ such that
\begin{align*}
\dot\gamma(t)=e^{i\theta(t)}
\end{align*}
and $0\le t_0<t_1\le L$ such that $\gamma$ restricted to 
$[0,t_0]$, $[t_0,t_1]$, and $[t_1,L]$
is a graph over the $x_1$-axis. 
\end{lemma}
\begin{proof}
We define the tangent angle $\theta:\: [0,L] \to \rz$ of $\gamma$ as the unique
continuous map such that $\theta(0)=\theta_0$ and
\begin{align*}
\dot\gamma(t)=e^{i\theta(t)} \text{ for all } t\in [0,L].  
\end{align*}
Since the curvature of $\gamma$ is positive, the tangent angle $\theta$ is
strictly increasing.
We apply Hopf's rotation angle theorem \cite{MR707850,MR1556906} to the simple positive oriented curve
$\gamma\oplus [-a,a]$ and find that the rotation angle of $\gamma\oplus [-a,a]$
is exactly $2\pi$. Consequently,
\begin{align*}
2\pi = \theta(L)+ (2\pi-\theta_L),  
\end{align*}
such that $\theta(L)=\theta_L$.
The curve $\gamma$ fails to be a graph over the $x_1$-axis, 
if $\theta(t)$ crosses $\pi/2$ or $3\pi/2$.
Since $\theta$ is strictly increasing,  
this can happen at most two times in the interval $(0,L)$.
This yields the claim. 
\end{proof}

\begin{lemma}
\label{l:min_estimate}
Let $\gamma \in C^2([0,L],\rz^2)$ be a unit speed curve with positive
geodesic curvature connecting
$(a,0)$ and $(-a,b)$, such that
\begin{align*}
\dot\gamma(t)=e^{i\theta(t)},  
\end{align*}
for some strictly increasing function $\theta \in C^0([0,L],\rz)$
satisfying $\pi/2 \le \theta(0)<\pi$ and $\pi< \theta(L)\le 3\pi/2$.
Then
\begin{align*}
\min \{k_\gamma(t)\where t \in [0,L]\} \le a^{-1}.
\end{align*}
\end{lemma}
\begin{proof}
Consider the upper half of the ball centered at $(0,0)$ and radius $a$
\begin{align*}
B_a^+ &:= \{(x,y)\in \rz^2 \where |x|\le a,\, y\ge 0,\, x^2+y^2\le a^2\},\\
C_a^+ &:= \{(x,\sqrt{a^2-x^2})\in \rz^2\where |x|\le a\},
\end{align*}
and 
\begin{align*}
s_1 &:= \sup\{s \in \rz\where (0,s)+\gamma \cap B_a^+ \neq \emptyset\}
\end{align*}
Obviously, there holds $s_1 \ge \max\{0,-b\}$. 
If $s_1>\max\{0,-b\}$, then
$s_1+\gamma$ and $B_a^+$ intersect in a point $(s_1,0)+\gamma(t_0)$
with $t_0\in (0,L)$ and $s_1+\gamma$ lies above $B_a^+$.
From the maximum principle the curvature of $\gamma$ 
at $\gamma(t_0)$ is smaller than $a^{-1}$.
If $s_1=0$, then $\gamma$ lies above $B_a^+$ and $\theta(0)$ has to be
$\pi/2$, such that
the slope
of $\gamma$ and $C_a^+$ coincide at $(a,0)$.
Writing $\gamma$ and $C_a^+$ as graphs over the $x_2$-axis
the maximum principle shows that
the curvature of $\gamma$ at $(a,0)$
is smaller than $a^{-1}$.
If $s_0=-b>0$ then $\theta(L)=3\pi/2$ and as above
we deduce 
$k_\gamma(L)\le a^{-1}$.
\end{proof}

\begin{lemma}
\label{l:max_estimate}
Let $\gamma \in C^2([0,L],\rz^2)$ be a unit speed curve with positive
geodesic curvature connecting
$(a,0)$ and $(-a,b)$, such that
\begin{align*}
\dot\gamma(t)=e^{i\theta(t)},  
\end{align*}
for some strictly increasing function $\theta \in C^0([0,L],\rz)$
satisfying $\pi/2 =\theta(0)$.
Moreover, if $b>0$, we assume that $\theta(L)=3\pi/2$, and if
$b \le 0$, we assume that $\pi< \theta(L) \le 3\pi/2$.
Then
\begin{align*}
\max \{k_\gamma(t)\where t \in [0,L]\} \ge a^{-1}.
\end{align*}
\end{lemma}
\begin{proof}
The curve $\gamma$ may be written as a graph over the interval
$[-a,a]$ for some function $g \in C^0([-a,a],\rz) \cap C^2((-a,a),\rz)$.
Let $G$ be set defined by
\begin{align*}
G:=\{(x,y)\in \rz^2\where -a\le x\le a,\, y\le g(x)\}.   
\end{align*}
Due to the positive curvature of $\gamma$ the set $G$ is convex and
\begin{align*}
G\cap \{(x,y)\in \rz^2\where x \in \{\pm a\},\, y> g(x)\}=\emptyset.  
\end{align*}
As in the proof of Lemma \ref{l:min_estimate} we consider
$C_a^+$ and
\begin{align*}
s_0 &:= \sup\{s \in \rz\where (0,s)+C_a^+ \cap G \neq \emptyset\},
\end{align*}
Obviously, there holds $s_0 \ge \max\{0,b\}$. 
If $s_0>\max\{0,b\}$, then
$(0,s_0)+C_a^+$ and $G$ intersect in a point $(t_0,g(t_0))$
with $|t_0|<a$ and $(0,s_0)+C_a^+$ lies above $G$.
From the maximum principle the curvature of $\gamma$ 
at $(t_0,g(t_0))$ is bigger than $a^{-1}$.
If $s_0=0$, then $C_a^+$ lies above $G$.
Since $\theta(0)=\pi/2$ the slope
of $\gamma$ and $C_a^+$ coincide at $(a,0)$.
From the maximum principle we deduce that
$k_\gamma(0)\ge a^{-1}$.
If $s_0=b>0$ then the slope
of $\gamma$ and $(0,b)+C_a^+$ coincide at $(-a,b)$
and the maximum principle shows that
$k_\gamma(L)\ge a^{-1}$.
\end{proof}

\begin{lemma}
\label{l:nonex_big}
Let $\gamma \in C^2([0,L],\rz^2)$ be a unit speed curve with positive
geodesic curvature connecting
$(a,0)$ and $(-a,0)$, such that the closed curve
$\gamma\oplus[-a,a]$
is simple and $\dot \gamma(L) \in \{e^{i\theta}\where \pi<\theta\le 5/2 \pi\}$.
If $\dot\gamma(0)=e^{-i\pi/2}$,
then the maximum $k_{max}$ and the minimum $k_{min}$
of the geodesic curvature of $\gamma$ satisfy
\begin{align*}
k_{min} \le \frac{k_{max}}{k_{max}a+1}.  
\end{align*}
\end{lemma}
\begin{proof}
We apply Lemma \ref{l:graph_small}, write
\begin{align*}
\dot\gamma(t)=e^{i\theta(t)}, \, -\pi/2< \theta(t)\le 5/2\pi,  
\end{align*}
and denote by $t_0$ the point such that
\begin{align*}
t_0 :=\sup\{t \in [0,L]\where \theta(s) \le \pi/2 \text{ for all }0\le s\le t\}.  
\end{align*}
By Lemma \ref{l:graph_small}
there holds $t_0<L$, $\theta(t_0)=\pi/2$, and $\theta(\cdot)$
is strictly increasing. Consequently, after a rotation by $\pi$,
we may apply Lemma \ref{l:max_estimate} and deduce that
$\gamma(t_0)=(x_0,y_0)$ with $x_0\ge a+2k_{max}^{-1}$.\\
We denote by $t_1$ the point
\begin{align*}
t_1 :=\sup\{t \in [t_0,L]\where \theta(s) < 3/2\pi \text{ for all }t_0\le s\le t\}.    
\end{align*}
Since $\gamma\oplus[a,-a]$ is simple, we have
$\gamma(t_1)=(x_1,y_1)$ for some $x_1\le -a$ (
if $t_1<L$, then $x_1<-a$). 
From Lemma \ref{l:min_estimate} applied to
$\gamma$ restricted to $[t_0,t_1]$ we see that
\begin{align*}
k_{min}\le \big(a+ k_{max}^{-1}\big)^{-1},  
\end{align*}
which yields the claim.
\end{proof}
We define for a given
curvature function $k\in C(\rz^2\times [0,1],\rz)$
the set of small solutions $L_{small}(k)$ and large solutions
$L_{large}(k)$ by 
\begin{align*}
L_{small}(k) &:= 
\{\gamma \in M_{small} \where \gamma \text{ solves \eqref{eq:1}.}\},\\
L_{large}(k) &:= \{\gamma \in M_{large} \where \gamma \text{ solves \eqref{eq:1}.}\}.
\end{align*}

\begin{lemma}
\label{l:set-compact}
Let $\{k_s \in C(\rz^2\times \rz,\rz^+)\where s \in [0,1]\}$ be a 
continuous family of prescribed curvature function, such that
\begin{align*}
\sup\{k_s(x,t)\where (x,t,s) \in \rz^2\times [0,1]^2\}a<1,\\
\inf\{k_s(x,t)\where (x,t,s) \in \rz^2\times [0,1]^2\}>0
\end{align*}
Then the set 
\begin{align*}
L_{small}:= 
\{\gamma \in M_{small} \where \gamma \text{ solves \eqref{eq:1} for some }k \in \{k_s\}\}
\end{align*}
is compact in $C^2([0,1],\rz^2)$. If, moreover, for all $s\in [0,1]$
\begin{align*}
\frac
{\sup_{(x,t) \in \rz^2\times [0,1]}\{k_s(x,t)\}}
{\sup_{(x,t) \in \rz^2\times [0,1]}\{k_s(x,t)\}a+1}
< 
\inf_{(x,t) \in \rz^2\times [0,1]}\{k_s(x,t)\}
\end{align*}
then
\begin{align*}
L_{large}:= \{\gamma \in M_{large} \where \gamma \text{ solves \eqref{eq:1} for some }k \in \{k_s\}\}
\end{align*}
is compact in $C^2([0,1],\rz^2)$.
\end{lemma}
\begin{proof}
We first show that that $L_{large}$ and $L_{small}$ are closed.
To this end we observe that any $\gamma \in L_{large} \cup L_{small}$
is parameterized proportional to its arclength.\\ 
Let $(\gamma_n)$ be a sequence in $L_{small}$ converging to $\gamma_0$
in $C^2([0,1],\rz^2)$. Choosing a subsequence, we may assume that
$\gamma_n$ is a solution to \eqref{eq:1} with $k=k_{s_n}$ for some
sequence $(s_n)$ converging to $s_0\in [0,1]$. Thus, $\gamma_0$
solves \eqref{eq:1} with $k=k_{s_0}$.
Using the maximum principle and the positive curvature of $\gamma_0$
it is easy to see that the curve $\gamma_0$ cannot touch itself
or the straight line $[-a,a]$ tangentially, such that
$\gamma_0\oplus [-a,a]$ remains simple as a limit of simple curves and
\begin{align*}
|\dot\gamma_0(0)|^{-1}\dot\gamma_0(0) \in \{e^{i \theta}\where 1/2 \pi \le \theta < \pi\},\\
\dot\gamma_0(1)|^{-1}\dot\gamma_0(1) \in \{e^{i \theta}\where \pi < \theta \le 3/2\pi\}.
\end{align*}
Since
\begin{align*}
\sup\{k_{s_0}(x,t)\where (x,t) \in \rz^3\}a<1
\end{align*}
by Lemma \ref{l:max_estimate} it is impossible that
\begin{align*}
|\dot\gamma_0(0)|^{-1}\dot\gamma_0(0)=e^{i \pi/2}
\text{ or }
|\dot\gamma_0(1)|^{-1}\dot\gamma_0(1)=e^{i 3\pi/2}.  
\end{align*}
Consequently, $\gamma_0$ is contained in $L_{small}$.\\
Let $(\gamma_n)$ be a sequence in $L_{large}$ converging to $\gamma_0$
in $C^2([0,1],\rz^2)$. As above, we may deduce that
$\gamma_0$ is a solution to \eqref{eq:1} with $k=k_{s_0}$ for some
$s_0\in [0,1]$, $\gamma_0\oplus [-a,a]$ is simple, and satisfies
\begin{align*}
|\dot\gamma_0(0)|^{-1}\dot\gamma_0(0) \in \{e^{i \theta}\where -\pi/2 \le \theta < \pi\},\\
\dot\gamma_0(1)|^{-1}\dot\gamma_0(1) \in \{e^{i \theta}\where \pi < \theta \le 5/2\pi\},
\end{align*}
and at least one of the following two conditions holds
\begin{align*}
{\dot\gamma(0)}{|\dot\gamma(0)|^{-1}} \in
\{e^{i\theta}\where -\pi/2\le \theta \le \pi/2\},\\
{\dot\gamma(1)}{|\dot\gamma(1)|^{-1}}  \in
\{e^{i\theta}\where 3/2\pi\le \theta \le 5\pi/2\}
\end{align*}
Using Lemma \ref{l:nonex_big} and the fact that
\begin{align*}
\frac
{\sup_{(x,t) \in \rz^2\times [0,1]}\{k_{s_0}(x,t)\}}
{\sup_{(x,t) \in \rz^2\times [0,1]}\{k_{s_0}(x,t)\}a+1}
< 
\inf_{(x,t) \in \rz^2\times [0,1]}\{k_{s_0}(x,t)\}  
\end{align*}
we exclude the possibility that
\begin{align*}
{\dot\gamma(0)}{|\dot\gamma(0)|^{-1}}=e^{-i\pi/2} \text{ or }
{\dot\gamma(1)}{|\dot\gamma(1)|^{-1}}=e^{i5\pi/2}.  
\end{align*}
If neither 
\begin{align*}
{\dot\gamma(0)}{|\dot\gamma(0)|^{-1}} \in
\{e^{i\theta}\where \theta < \pi/2\}  
\end{align*}
nor
\begin{align*}
{\dot\gamma(1)}{|\dot\gamma(1)|^{-1}}  \in
\{e^{i\theta}\where 3/2\pi < \theta\}
\end{align*}
then Lemma \ref{l:max_estimate} leads to a contradiction.
Thus, $\gamma_0$ belongs to $L_{large}$.\\
To show the compactness of $L_{large}$ and $L_{small}$ we fix a 
sequence $(\gamma_n)$ of solutions in $L_{large} \cup L_{small}$.
Since $\gamma_n \oplus [-a,a]$ is simple we may apply the Gau\ss-Bonnet formula
and get
\begin{align*}
2\pi = \alpha_{1,n}+\alpha_{2,n}+\int_{\gamma_n} k_{\gamma_n},  
\end{align*}
where $\alpha_{1,n},\, \alpha_{2,n} \in (-\pi/2,\pi)$ are the outward angles
at $t=0$ and $t=1$ of the piecewise $C^2$ curve $\gamma_n \oplus [-a,a]$.
Consequently,
\begin{align*}
L(\gamma_n) \inf_{(x,t,s) \in \rz^2\times[0,1]^2}\{k_s(x,t)\}
\le \int_{\gamma_n} k_{\gamma_n} \le 3\pi,  
\end{align*}
where $L(\gamma_n)$ denotes the length of $\gamma_n$. Hence,
$L(\gamma_n)$ is uniformly bounded, which yields a uniform
$C^1$-bound of $\gamma_n$. Using the equation \eqref{eq:1} 
and the Arzela-Ascoli theorem we may extract
a subsequence of $(\gamma_n)$, which converges in $C^2([0,1],\rz^2)$.
This finishes the proof.    
\end{proof}

\section{The Leray-Schauder degree}
\label{sec:degree}
For $a>0$ we consider the affine space
\begin{align*}
C^2_{a,-a}([0,1],\rz^2):=\Big\{\gamma \in C^2([0,1],\rz^2):\: 
\gamma(0)=\begin{pmatrix} a\\ 0\end{pmatrix}
\text{ and }
\gamma(1)=\begin{pmatrix} -a\\ 0\end{pmatrix}
\Big\}.   
\end{align*}
The operator $L_k$ is defined by
\begin{align*}
L_k:\: C^2_{-a,a}([0,1],\rz^2) \to C^2_{-a,a}([0,1],\rz^2)
\end{align*}
\begin{align*}
L_k(\gamma) := \big(-D_t^2\big)^{-1} \Big(-\ddot \gamma 
+ |\dot \gamma(\cdot)| k(\gamma(\cdot),\cdot) J(\dot\gamma(\cdot))\Big),  
\end{align*}
where the operator $D_t^2$ is considered as an isomorphism
\begin{align*}
D_t^2:\: C^2_{-a,a}([0,1],\rz^2) \to C^0([0,1],\rz^2).  
\end{align*}
Since
\begin{align*}
|\dot \gamma(\cdot)| k(\gamma(\cdot),\cdot) J(\dot\gamma(\cdot))\in C^0([0,1],\rz^2)   
\end{align*}
depends only on $\gamma$ and $\dot\gamma$, the map 
\begin{align*}
\gamma \mapsto   
\big(-D_t^2\big)^{-1} \Big(|\dot \gamma(\cdot)| k(\gamma(\cdot),\cdot) J(\dot\gamma(\cdot))\Big)
\end{align*}
is compact from $C^2_{-a,a}([0,1],\rz^2)$ to itself. Thus $L_k$ is of the
form $Id - \text{compact}$ and the Leray-Schauder degree of $L_k$ is defined.\\
Fix $a>0$ and a function $k \in C(\rz^2\times [0,1],\rz)$ satisfying
\begin{align*}
0<\inf_{\rz^2\times [0,1]}k \le \sup_{\rz^2\times [0,1]} k <a^{-1},\\
\frac
{\sup_{(x,t) \in \rz^2\times [0,1]} k(x,t)}
{\sup_{(x,t) \in \rz^2 \times [0,1]} k(x,t)a+1}
< 
\inf_{(x,t) \in \rz^2\times [0,1]} k(x,t).  
\end{align*}
We define for $s \in [0,1]$ the function $k_s \in C^0(\rz^2\times [0,1],\rz)$ by
\begin{align*}
k_s(x,t):= (1-s)\big(\sup_{(x,t) \in \rz^3} k(x,t)\big)+s k(x,t).  
\end{align*}
Then the family $\{k_s\where s\in [0,1]\}$
satisfies the assumptions of Lemma \ref{l:set-compact} and the sets $L_{large}$ and
$L_{small}$ are compact. Thus, there is $R>0$ such that
\begin{align*}
L_{large}\cup L_{small} \subset \{\lambda \in C^2([0,1],\rz^2)
\where \|\lambda\|_{C^2([0,1],\rz^2)}<R\}.  
\end{align*}
Consequently, if we define the open sets
\begin{align*}
M_{small,R}:= \{\lambda \in M_{small} \where \|\lambda\|_{C^2([0,1],\rz^2)}<R\},\\
M_{large,R}:= \{\lambda \in M_{large} \where \|\lambda\|_{C^2([0,1],\rz^2)}<R\},
\end{align*}
then from the homotopy invariance of the degree
\begin{align}
\label{eq:5}
\deg(L_{k},M_{small,R},0)&= \deg(L_{k_0},M_{small,R},0), \notag\\
\deg(L_{k},M_{large,R},0)&= \deg(L_{k_0},M_{large,R},0).  
\end{align}
To compute the degree of $L_{k_0}$ we note that solutions to \eqref{eq:1}
with a constant function $k_0$ are given by curves with constant geodesic curvature
$k_0$, i.e. subarcs of a $n$-fold iterate of a circle with radius $k_0^{-1}$. 
Thus the required simplicity and the bounds on the slope yields
\begin{align*}
L_{small}(k_0) &= 
\{\gamma_s(t):= k_0^{-1} e^{i(\alpha_0+\omega_s t)}-i k_0^{-1}\sin(\alpha_0)\},\\
L_{large}(k_0) &= 
\{\gamma_b(t):= k_0^{-1} e^{i(-\alpha_0+\omega_b t)}+i k_0^{-1}\sin(\alpha_0)\},  
\end{align*}
where
\begin{align*}
\alpha_0 := \arccos(k_0 a) \in (0,\pi/2),\\
\omega_s := \pi-2\alpha_0 \in (0,\pi),\\
\omega_b := \pi+2\alpha_0 \in (\pi,2\pi).  
\end{align*}
Consequently, we have
\begin{align}
\label{eq:7}
\deg(L_{k},M_{small,R},0)&= \deg_{loc}(DL_{k_0}\eval_{\gamma_{s}},0), \notag\\
\deg(L_{k},M_{large,R},0)&= \deg_{loc}(DL_{k_0}\eval_{\gamma_{b}},0).   
\end{align}
To compute the local degree's we note
for $V \in C_{0,0}^{2}([0,1], \rz^2)$ and $*\in \{s,b\}$
\begin{align*}
DL_{k_0}\eval_{\gamma_{*}}(V) &=(-D_t^2)^{-1}
\big(-\ddot V
+\langle \dot\gamma_{*},\dot V\rangle |\dot\gamma_{*}|^{-1} k_0 J(\dot \gamma_{*})
+|\dot\gamma_{*}| k_0 J(\dot V)\big)\\
&= 
(-D_t^2)^{-1}
\big(-\ddot V
-\omega_* \langle i e^{i(\alpha_0+\omega_* t)},\dot V\rangle e^{i(\alpha_0+\omega_* t)}\\
&\qquad +\omega_* J(\dot V)\big).   
\end{align*}
For $\lambda \in [-1,1]$ we consider the family of operators 
$A_\lambda: C_{0,0}^{2}([0,1], \rz^2) \to C_{0,0}^{2}([0,1], \rz^2)$ defined by
\begin{align*}
A_\lambda (V) := (-D_t^2)^{-1}&
\big(-\ddot V
-(1-\lambda)
\omega_* \langle i e^{i(\alpha_0+\omega_* t)},\dot V\rangle e^{i(\alpha_0+\omega_* t)}\\
&\qquad +(1+\lambda) \omega_* J(\dot V)\big).
\end{align*}
Writing
\begin{align}
\label{eq:4}
V(t) = \alpha(t) e^{i(\alpha_0+\omega_* t)} + \beta(t) i e^{i(\alpha_0+\omega_* t)},
\end{align}
for some $\alpha,\, \beta \in C_0^2([0,1],\rz)$ we find
\begin{align*}
A_{\lambda,*}(V) 
=(-D_t^2)^{-1}
\Big(
&\big(
-\ddot \alpha(t)- \omega_*^2 \alpha(t)
\big)
e^{i(\alpha_0+\omega_* t)}\\
&+
\big(
-\ddot\beta(t)-(1-\lambda) \omega_*\dot \alpha(t)- \lambda \omega_*^2\beta(t) 
\big)
ie^{i(\alpha_0+\omega_* t)}
\Big)\\
\end{align*}
The eigenvalues of the problem
\begin{align*}
\ddot \phi(t) = \lambda \phi(t) \text{ for }t \in [0,1] \text{ and }\phi(0)=\phi(1)=0  
\end{align*}
are given by 
\begin{align}
\label{eq:3}
\{\pi^2n^2\where n\in \nz\}.  
\end{align}
Since $\omega_s < \pi$, each $A_{\lambda,s}$ is injective
and due to its form, identity-compact, $A_{\lambda,s}$ is invertible
for each $\lambda \in [0,1]$. By the homotopy invariance of the degree we obtain
\begin{align}
\label{eq:6}
\deg_{loc}(DL_{k_0}\eval_{\gamma_{s}},0)=
\deg_{loc}(A_{1,s},0) = \deg_{loc}(id,0)=1,
\end{align}
where we used for the second equality 
the admissible homotopy $\{B_{\sigma}\where \sigma \in [0,1]$ given by
\begin{align*}
B_{\sigma}(V) := (-D_t^2)^{-1}
\big(-\ddot V +2(1-\sigma)\omega_* J(\dot V)\big).  
\end{align*}
To compute the degree of $D L_{k_0}\eval_{\gamma_b}$ we note that by the above analysis
and the homotopy property we may replace $k_0$ by some constant $k_1$ close to 
$a$ without changing the degree, such that we may assume
\begin{align}
\label{eq:2}
\pi <\omega_b< \sqrt{2}\pi.   
\end{align}
Moreover, using the homotopy $\{A_{\lambda,b}\where \lambda \in [-1,0]\}$, we
see that 
\begin{align*}
\deg_{loc}(D L_{k_0}\eval_{\gamma_b},0) =
\deg_{loc}(A_{-1,b},0).   
\end{align*}
To compute $\deg_{loc}(A_{-1,b},0)$ we consider the decomposition
\begin{align*}
C^2_{0,0}([0,1],\rz^2)= U_1 \oplus U_2,   
\end{align*}
where
\begin{align*}
U_1 &:= \{V \in C^2_{0,0}([0,1],\rz^2) \where 
\int_0^1 V(t) \cdot \big(\sin(\pi t)e^{i(\alpha_0+\omega_b t)}\big) \, dt =0\},\\
U_2 &:= \text{span}(\sin(\pi t)e^{i(\alpha_0+\omega_b t)}).  
\end{align*}
Using the decomposition in \eqref{eq:4} 
we fix $V_1 \in U_1\setminus \{0\}$ and $V_2 \in U_2\setminus \{0\}$,
\begin{align*}
V_1(t) &= \alpha(t) e^{i(\alpha_0+\omega_b t)} + \beta(t) i e^{i(\alpha_0+\omega_b t)},\\
V_2(t) &= \lambda \sin(\pi t)e^{i(\alpha_0+\omega_b t)}.  
\end{align*}
From \eqref{eq:3} and \eqref{eq:2} we obtain
\begin{align*}
\langle D_t &A_{-1,b}(V_1),D_t V_1\rangle_{L^2([0,1],\rz^2)}\\
&= \langle -(D_t)^2 A_{-1,b}(V_1),V_1\rangle_{L^2([0,1],\rz^2)}\\
&= \int_0^1 \big(-\ddot \alpha(t)-\omega_b^2\alpha(t)\big)\alpha(t)+
\big(-\ddot \beta(t)-2\omega_b\dot\alpha(t)+ \omega_b^2\beta(t)\big)\beta(t)\, dt\\
&= \int_0^1 (\dot\alpha(t))^2-2\omega_b^2(\alpha(t))^2
+(\dot \beta(t)-\omega_b\alpha(t))^2+\omega_b^2 (\beta(t))^2 \,dt\\
&\ge {(4\pi^2-2\omega_b^2)}(\alpha(t))^2 + \omega_b^2 (\beta(t))^2\\
&> 0, 
\end{align*}
\begin{align*}
\langle D_t A_{-1,b}(V_2),D_t V_2\rangle_{L^2([0,1],\rz^2)}
&= \lambda^2\int_0^1 (\pi^2-\omega_b^2)(\sin(\pi t))^2\, dt\\
&= \frac12 \lambda^2 (\pi^2-\omega_b^2)<0,  
\end{align*}
and
\begin{align*}
\langle D_t A_{-1,b}(V_1),D_t V_2\rangle_{L^2([0,1],\rz^2)}
&= \lambda \int_0^1 \big(-\ddot \alpha(t)-\omega_b^2\alpha(t)\big) \sin(\pi t) dt\\
&= \lambda (\pi^2-\omega_b^2) \int_0^1 \alpha(t) \sin(\pi t) dt=0.  
\end{align*}
Thus the following homotopy is admissible
\begin{align*}
[0,1]\ni \sigma \mapsto \sigma C + (1-\sigma) A_{-1,b},  
\end{align*}
where $C \in \mathcal{L}(C^2_{0,0}([0,1],\rz^2),C^2_{0,0}([0,1],\rz^2))$ is given
in the decomposition $U_1 \oplus U_2$ by
\begin{align*}
C:=
\begin{pmatrix}
id &0\\
0 &-1  
\end{pmatrix}
.
\end{align*}
From the above computations we finally see that
\begin{align*}
\deg_{loc}(D L_{k_0}\eval_{\gamma_b},0) 
= \deg_{loc}(C,0)=-1,   
\end{align*}
which yields together with \eqref{eq:7} the proof of Theorem \ref{intro:main_theorem}
announced in the introduction.

\bibliographystyle{plain}
\bibliography{planar_plateau}

\end{document}